\documentclass{elsart}

\usepackage{amssymb}
\usepackage{bm}
\usepackage{epsfig}
\usepackage{amsfonts}
\usepackage{color}
\begin{document}

\begin{frontmatter}



\title{Multiscale analysis of re-entrant production lines: An equation-free approach}


\author{Y. Zou\corauthref{a}}
\author{$\quad$I.G. Kevrekidis\corauthref{b}}
\author{$\quad$D. Armbruster\corauthref{c}}
\address[a]{Department of Chemical Engineering and PACM, Princeton University, Princeton, NJ 08544, USA, (yzou@Princeton.EDU).}
\address[b]{Department of Chemical Engineering and PACM, Princeton University, Princeton, NJ 08544, USA, (yannis@Princeton.EDU). To whom correspondence should be addressed.}
\address[c]{Department of Mathematics, Arizona State University, Tempe, AZ 85287-1804,
USA, (armbruster@asu.edu).}
\begin{abstract}
The computer-assisted modeling of re-entrant production lines, and, in particular, simulation scalability,
is attracting a lot of attention due to the importance of such lines in semiconductor manufacturing.
Re-entrant flows lead to competition for processing capacity among the items produced,
which significantly impacts their throughput time (TPT).
Such production models naturally exhibit two time scales: a short one, characteristic
of single items processed through individual machines, and a longer one,
characteristic of the response time of the entire factory.
Coarse-grained partial differential equations for the spatio-temporal evolution of
a ``phase density" were obtained through a kinetic theory approach in Armbruster {\it et al.} \cite{Armbruster:04}.
We take advantage of the time scale separation to directly solve such coarse-grained
equations, even when we cannot derive them explicitly, through an equation-free
computational approach.
Short bursts of appropriately initialized stochastic fine-scale simulation are used
to perform coarse projective integration on the phase density.
The key step in this process is {\it lifting}: the construction of fine-scale, discrete
realizations consistent with a given coarse-grained phase density field.
We achieve this through computational evaluation of conditional distributions of
a ``phase velocity" at the limit of large item influxes.

\end{abstract}

\end{frontmatter}

\bigskip
\bigskip
\bigskip
\bigskip

{\bf 1$\quad$Introduction}

Semiconductor production lines are billion dollar investments 
and their performance characteristics are carefully studied and
modelled in detail through discrete event simulations.
In such simulations, a network of processors (machines) is set up with connections
that describe the product flow through a factory.
As the production process of a typical semiconductor goes through several similar layers,
machines may process an item several times at different stages in the production process (re-entrancy).
The passage of an item through a factory is characterized by its {\it throughput time} (TPT),
which may vary due to several factors, most importantly
the stochastic features of the machines (their failure statistics)
and the interaction with human operators.
For a typical discrete event model, therefore, the time through a factory process is
a random variable sampled from a (given) probability distribution.

A basic model describing the flow of a single type of products through a processing factory
is given by
\begin{eqnarray}
   (a) &\quad& e_n = a_n + \tau_n, \nonumber \\
   (b) &\quad& d {\mathcal P} \{ \tau_n \le r \} = {\mathcal T}(r,t=a_n) dr.
\label{simplemodel:eqn}
\end{eqnarray}
where $a_n$ is the time when item $n$ enters the factory and $e_n$ its exit time.

The time interval $\tau_n$ for which
the item stays in the processing factory is sampled from the distribution ${\mathcal T}(r,t=a_n)$,
which is determined at the entrance time.
In engineering practice, the most significant influence on the form of this distribution
comes from the total number of items in progress, i.e., the {\it work in progress} (WIP) 
(see \cite{Armbruster:02,Armbruster:04}).
As a result, one rewrites ${\mathcal T}(r,t=a_n)$ in the form of ${\mathcal T}(r,WIP)$,
where ${\it WIP}$ is the current WIP at the time $a_n$.
When we fix the TPT at the beginning of the simulation,
we essentially treat the factory as a single queue whose length at item arrival
determines the time that the item needs to get processed.
However, for re-entrant flows, significant changes
in WIP during the time interval  $[a_n,e_n]$ may lead to a change in the TPT.
In \cite{Armbruster:04}, this situation is treated by introducing
the concepts of {\it phase} (a scaled position) and {\it phase velocity}.
The latter is stochastically updated from a distribution that depends on the
total WIP in the factory.
In this manner, the TPT is a variable that can dynamically change,
impacted by later arrivals of new items.

In \cite{Armbruster:04} an advection-diffusion equation for the
phase density $\rho(x,t)$ as a function of the phase, $x$, and the
time $t$ is derived via a Chapman-Enskog expansion of a Boltzmann
equation.
This PDE can be solved in a straightforward way allowing
us to determine all relevant quantities: total Work in Progress (WIP),
the WIP distribution, the throughput and the throughput time.
The relevant time scale for the PDE simulation is the time scale of the TPT evolution
of the whole factory.
In contrast,
for a discrete event simulation model, the relevant time scale is that of processing by an
individual machine, which is at least an order of magnitude smaller than the TPT. In addition, in order to obtain a meaningful statistical average, a large ensemble of items must be used.
As a result, the discrete event simulations are time consuming, and the simulation effort
scales with the number of steps in the factory and the number of items produced. In contrast, the continuous model equation is deterministic and
its simulation is independent of the number of steps and the number of items produced.

Comparison of a discrete event simulation (DES) model detailed below and
the explicitly derived advection-diffusion equation was performed in \cite{Armbruster:04}
and the results were in reasonable agreement; in that case, the PDE model
was derived ``from first principles''.
However, the particular discrete event simulations
already constituted a significant approximation, as they were written for particles in a continuum.
In the general case of a large scale DES with many items and many machines,
global continuum evolution equations for the item density cannot easily
be analytically derived in a closed form.
The question then arises whether we can still observe the time evolution
at the density level, maintaining the scalability advantages of a density-type model.
The recently developed Equation-Free (EF) approach is aimed at precisely this type
of problem: we have a fine-scale (here DES) simulator, we suspect that we can model
its coarse-grained behavior at the level of a density evolution equation, but do
not have this equation explicitly available.
The approach attempts to solve the equation (integrate it, find its stationary solutions
and their stability, etc.) by designing short bursts of appropriately initialized simulations
with the fine-scale (here DES) model.
The quantities required for scientific computation with the unavailable density-level
equation are then {\it estimated on demand} from these short fine-scale runs
\cite{Theodoropoulos:00,Kevrekidis:03,Kevrekidis:04}.
In this paper we will use a particular equation-free algorithm (coarse projective
integration  \cite{Gear:01,GearA:01,Rico-Martinez:04}) to illustrate the application of the approach to the re-entrant line problem.

Other computational tasks like coarse-grained fixed point computation, stability, bifurcation
analysis, control, and computation of self-similar solutions, have been demonstrated in the
literature for various types of ``inner", fine-scale simulators \cite{Gear:02,Makeev:02,MakeevA:02,Chen:03,Zou:05}.
A detailed discussion of the methods can be found in \cite{Kevrekidis:03,Kevrekidis:04}.
A key step in EF computations is the so-called ${\it lifting}$: the
construction of fine-scale states consistent with given values of their coarse-scale
{\it observables} (the dependent variables in the unavailable coarse equation).
This is required in order to initiate new short bursts of fine-scale evolution.
This is not a deterministic step - many fine-scale configurations share the same observables.
However, the assumption that an evolution equation exists and meaningfully closes at the level of these
coarse observables suggests that the details of various lifting realizations
do not affect the long-time coarse-scale evolution of the observables
(see \cite{Gear:02} for detailed discussions).
In applying the EF methods to re-entrant factory production, one may use the joint
``phase and phase velocity" density $f(x,v,t)$ as the observable of choice.
However, the work in \cite{Armbruster:04} suggests that at an appropriate (large influx) parameter limit,
one can close an equation in terms of a simpler coarse-grained observable: just
the phase density $\rho(x,t)$.
While for the former the lifting procedure would be straightforward
(just generating two random variables according to their joint probability density function \cite{Zou:05}), for the latter care must be taken in generating joint probability densities based on just
the phase density $\rho(x,t)$.

The paper is organized as follows.
In Section 2 we briefly describe the phase model of a re-entrant factory discussed in \cite{Armbruster:04}.
Section 3 contains our main equation-free results: our lifting procedure, as well as
coarse projective integration based on ``finite-difference''-type observations
of the phase density evolution.
We conclude with a brief summary, discussion, and possible extensions of the
approach.

{\bf 2$\quad$A Phase Model for a Re-entrant Factory} \label{phasemodel:sec}

{\it 2.1$\quad$The Discrete Model} \label{discrete:sec}

The phase ${\it s}$ of an item in a given factory of a supply chain
is defined as the antiderivative of the phase velocity,
$1/\tau(t)$, where $\tau(t)$ may change with time,
- in contrast to the constant $\tau_n$ in Eqn. (\ref{simplemodel:eqn}).
With a deterministic TPT $\tau(t)$ the improved model reads:
\begin{eqnarray}
   (a) &\quad& s=\phi(t), \nonumber \\
   (b) &\quad& {{d\phi} \over {dt}} = { 1 \over \tau(t)}, \quad \phi(a_n)=0, \quad t \ge a_n.
\label{phasemodel:eqn}
\end{eqnarray}
The exit time $e_n$ is the time at which $\phi(t)=1$.
When $\tau(t)$ is constant, model (\ref{phasemodel:eqn}) reduces to model (\ref{simplemodel:eqn}).
When the TPT is a random variable sampled from a distribution ${\mathcal T}(r,t)$,
a discrete form corresponding to (\ref{phasemodel:eqn})
must be used (where $\omega$ is an update frequency):
\begin{eqnarray}
 &\quad& \phi(t+{1 \over \omega}) = \phi(t) + {1 \over {\omega \tau(t)}}, \nonumber \\
 &\quad& \phi(a_n) = 0, \quad t \ge a_n, \quad d {\mathcal P} \{ \tau(t) \le r \} = {\mathcal T}(r,t)dr.
\label{discrete1:eqn}
\end{eqnarray}
Qualitatively, the value of $\omega$ is influenced by the number of items entering the factory
within a characteristic TPT scale.

Given a prescribed TPT distribution function ${\mathcal T}(r,t)$ (which
we assume can be expressed in terms of WIP as ${\mathcal T}(r,WIP)$), an
executable algorithm for evolving the phase ${\it s}$ and the
real-time TPT $\tau(t)$ as well as the exit time $e_n$ is also given
in \cite{Armbruster:04}, namely,
\begin{eqnarray}
  (a) &\quad& \phi(t+\Delta t) = \phi(t) + { {\Delta t} \over {\tau(t)}}, \nonumber \\
      &\quad& \tau(t+\Delta t) = \kappa(t) \eta(t) + (1-\kappa(t)) \tau(t), \quad t \ge a_n, \nonumber \\
  (b) &\quad& {\mathcal P} \{ \kappa(t) = 1 \} = \omega \Delta t, \quad {\mathcal P} \{ \kappa(t) = 0 \} = 1- \omega \Delta t, \nonumber \\
      &\quad& d{\mathcal P} \{ \eta(t) \le r \} = {\mathcal T}(r,t)dr, \nonumber \\
  (c) &\quad& \phi(a_n)=0, \quad   d{\mathcal P} \{ \tau(a_n) \le r \} = {\mathcal T}(r,a_n)dr,
\label{discrete2:eqn}
\end{eqnarray}
where the update frequency $\omega$ depends on the throughput time
$\tau(t)$ and the time $t$ itself, i.e., $\omega =
\omega(\tau(t),t)$.
The simulation time step $\Delta t$ is chosen
such that $\Delta t < { 1 \over {\omega}}$.
For small enough $\Delta t$ this algorithm can be used to
numerically solve model (\ref{discrete1:eqn}).

In \cite{Armbruster:04}, $\omega$ is chosen as
\begin{equation}
   \omega(r,t)= { {\lambda(t) } \over {r T_{-1}}},
\label{omega:eqn}
\end{equation}
where $\lambda(t)$ is the influx of items and $T_{-1}= \int r^{-1} {\mathcal T}(r,t) dr$.

The procedure to compute the phase and throughput time (TPT) of an
item in the factory is then as follows:
\begin{enumerate}
\item Set the item's initial phase $\phi$ to 0 as it enters the factory; increase the WIP $W(t)$ by one and adjust the distribution of
TPT, ${\mathcal T}(r,t)$ accordingly (remember that we have assumed that
this distribution can be written as ${\mathcal T}(r,WIP)$); the item's
initial TPT is then sampled from the udpated ${\mathcal T}(r,WIP)$;
\item Compute $\omega(\tau(t),t)$ as above, and sample the parameter $\kappa(t)$
according to its distribution (\ref{discrete2:eqn})(b);
\item Compute the TPT $\tau(t+\Delta t)$ and the phase $\phi(t+\Delta t)$ according to (\ref{discrete2:eqn})(a);
\item Update the distribution of TPT, ${\mathcal T}(r,t)$, according to current WIP at the time $t+\Delta t$, 
and go back to Step (2); repeat this loop
until $\phi=1$. The time when  $\phi=1$ is the exit time of the item, $e_n$.
\end{enumerate}

{\it 2.2$\quad$Density Equations} \label{density:sec}

Let the joint number density of phase and phase velocity,
$f(x,r,t)$, be defined as $f(x,r,t)= {{\partial^2 F(\phi \le x, \tau
\le r, t)} \over {\partial \phi \partial r}}$, where $ F(\phi \le x,
\tau \le r, t)$ is the number of items whose phase $\phi \le x$ and
TPT $\tau \le r$.
The PDE for $f(x,r,t)$ is derived in \cite{Armbruster:04}.
It is given by
\begin{eqnarray}
&\quad&  {{\partial f} \over {\partial t}} + {1 \over r}  {{\partial f} \over {\partial x}} = {\mathcal T}(r,t) \int \omega(r',t)f(x,r',t)dr' - \omega(r,t)f(x,r,t), \nonumber \\
&\quad& x>0, \quad t>0 \nonumber \\
&\quad&  f(0,r,t)= r \lambda(t) {\mathcal T}(r,t), \quad f(x,r,0)=0
\label{jointdensity:eqn}
\end{eqnarray}
weakly in $x$, $r$ and $t$, where $\lambda(t)$ is the influx of items.
The chosen form of $\omega$ in (\ref{omega:eqn}) guarantees that the
discrete model (\ref{simplemodel:eqn}), when there are no temporal changes in influx and TPT,
 is a Monte Carlo scheme
corresponding to a particular solution of (\ref{jointdensity:eqn}), and that the
instances of random number generation are roughly the same for
models (\ref{simplemodel:eqn}) and (\ref{discrete2:eqn}).

If we define the number density of phase, $\rho(x,t)$, as
$\rho(x,t)= \int f(x,r,t)dr$, then in the limit that $\lambda_0
\sigma_{\mathcal T}^0 \gg 1$ (where $\lambda_0$ and $\sigma_{\mathcal T}^0$
are characteristic scales of influx and standard deviation of TPT,
respectively) one can write a closed equation for its evolution as
follows:
\begin{eqnarray}
&\quad& { {\partial \rho} \over {\partial t} } +  { {\partial F} \over {\partial x} } = 0, \quad x,t>0, \quad F(x,t)=C(t)\rho - D(t){ {\partial \rho} \over {\partial x} }, \nonumber \\
&\quad& C= {1 \over T_1} + {T_{-1} \over \lambda} {1 \over T_1} { {\partial (T_2/T_1)} \over {\partial t}}, \quad D= {T_{-1} \over \lambda} {{T_2 - T_1^2} \over T_1^3}, \nonumber \\
&\quad& F(0,t)=\lambda(t), \quad \rho(x,0)=0
\label{density:eqn}
\end{eqnarray}

The moments of the probability distribution ${\mathcal T}(r,t)$
appearing in this formula, $T_i(t),i=-1,1,2$, are defined by $T_i(t)
= \int r^i {\mathcal T}(r,t) dr, i=-1,1,2$.

Based on the proof of Theorem T3 in \cite{Armbruster:04}, it can
also be deduced that in this limit, the joint number density
$f(x,r,t)$ is controlled by (slaved to) just the number density of
phase, $\rho(x,t)$.
Their relationship is given by
\begin{equation}
   f(x,r,t) = { {\rho r {\mathcal T}(r,t)} \over T_1}
\label{reducedrelation:eqn}
\end{equation}

Equation (\ref{reducedrelation:eqn}) gives a simple relationship
between the joint number density and the number density of phase.
It implies that the conditional number density $f(r|x,t)$ ($= f(x,r,t)/ \rho(x,t)$)
is a constant $ {r {\mathcal T}(r,t)}/ {T_1}$ with respect to the phase coordinate $x$.

{\bf 3$\quad$Equation-Free Analysis for the Two-scale System} \label{EFanalysis:sec}

The Equation-Free approach consists of an ensemble of computational
tools that can be used to study the coarse-grained, macroscopic behavior
of systems based on their underlying fine-scale, microscopic
simulators (see e.g.
\cite{Gear:01,Gear:02,Makeev:02,Makeev:02,Chen:03,Zou:05}).
The basic element of these algorithms is the so-called {\it coarse
time-stepper}, whose role is to connect observables across different
scales.
The coarse time-stepper consists essentially of three components:
{\it lifting} and {\it micro-simulation} followed by {\it
restriction}.
The {\it lifting} is a procedure that generates micro-scale
realizations of a system state  consistent with given values of
their macro-scale observables; while {\it restriction} is the
reverse: obtaining macro-scale observables from the fine-scale
system state.
Since a fine-scale system state normally possesses far more degrees
of freedom than a few macro-scale observables, the lifting procedure
is not a one-to-one mapping in general.
Care needs to be taken when a lifting algorithm is implemented and
tests are usually required to check if the macro-scale evolution is
sensitive to the details of a particular lifting.

If we have reason to believe that useful macroscopic equations
accurately {\it close} at the level of a few macroscopic observables
(e.g. in terms of a few moments of microscopically evolving
distributions, as is the case in hydrodynamic equation derivation
from Boltzmann-level descriptions), equation-free schemes allow us
to solve the coarse-grained equations without the explicit closures
needed to write them in closed form.
The idea is to use short bursts of fine-scale simulation to evaluate
the right hand sides of the unavailable closed equations {\it on
demand}.
Traditional continuum numerical analysis is thus transformed into
protocols for the design and processing of repeated short bursts of
computational experiments with the fine-scale solver.
The assumption that coarse-grained equations in principle exist and
close in terms of a few coarse-grained observables implicitly
suggests that the details of a lifting should be quickly forgotten;
the role of the brief fine-scale simulation is to ``implement'' this
loss of memory of initial features of the fine-scale state.

{\it 3.1$\quad$The Lifting Step in the Equation-Free Approach} \label{lifting:sec}

We are interested in enabling short bursts of simulation with the
discrete phase model (\ref{discrete2:eqn}) to numerically analyze
the evolution of the number density $\rho(x,t)$ - a coarser-grained
observable than the full $f(x,r,t)$.
This should only be attempted at conditions when a deterministic
equation closes with $\rho(x,t)$, i.e. when the evolution of the
full $f(x,r,t)$ is controlled by $\rho(x,t)$.
For such a closed equation to exist, it should be possible (possibly
after a short initial transient) to express $f(x,r,t)$ in terms of
$\rho(x,t)$. This section investigates the effect of the parameter $\lambda_0
\sigma_{\mathcal T}^0$ on the relationship between $\rho(x,t)$ and
$f(x,r,t)$.
This is both because $\lambda_0 \sigma_{\mathcal T}^0$ is the only
dimensionless parameter that appears in the discrete model, and also
because we are interested in cases of large item influx within the
characteristic TPT scale, when $\lambda_0 \sigma_{\mathcal T}^0 \gg 1$.

\begin{figure}
\begin{minipage}{4.5cm}
\epsfig{figure=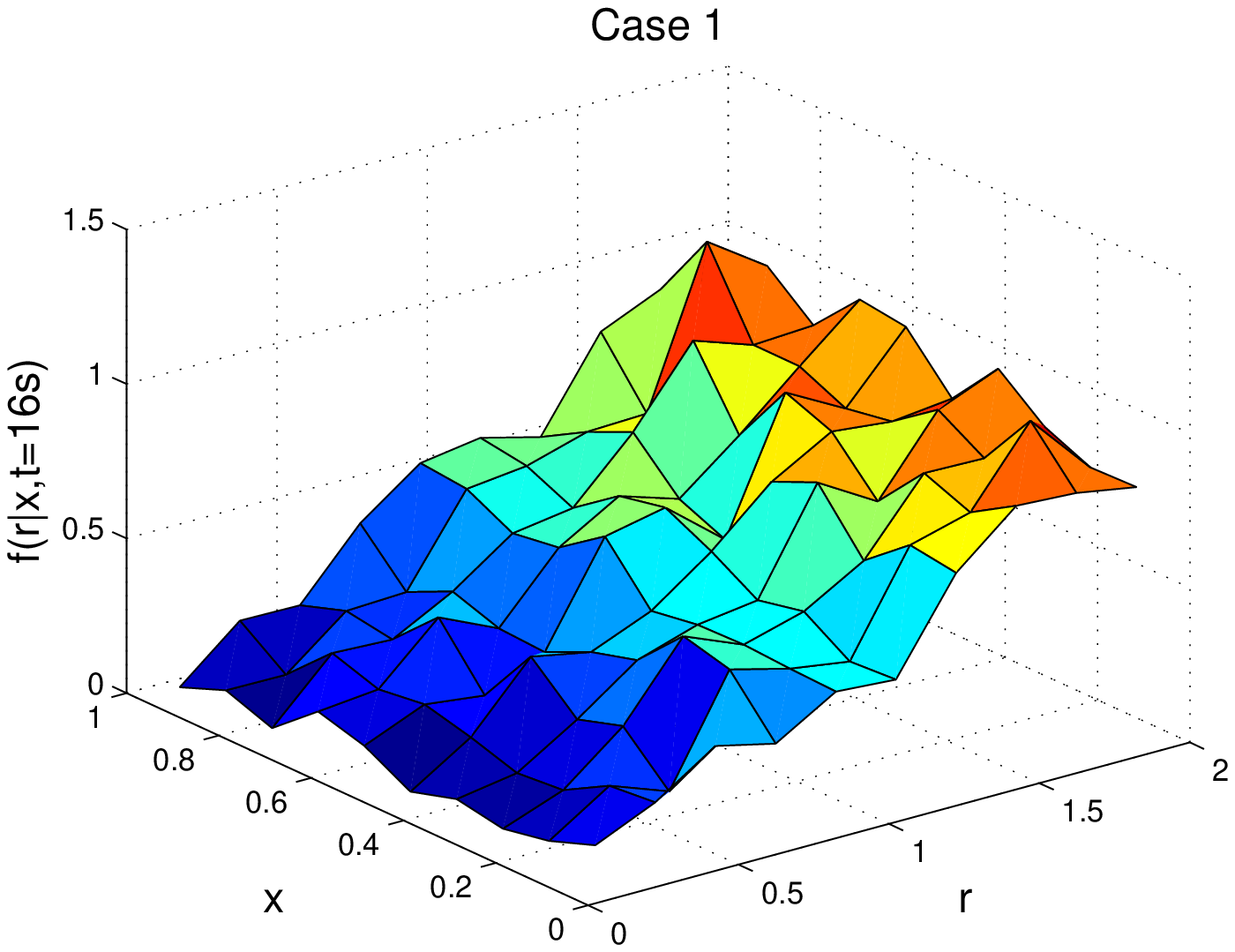,width=\textwidth}
\end{minipage}
\hfill
\begin{minipage}{4.5cm}
\epsfig{figure=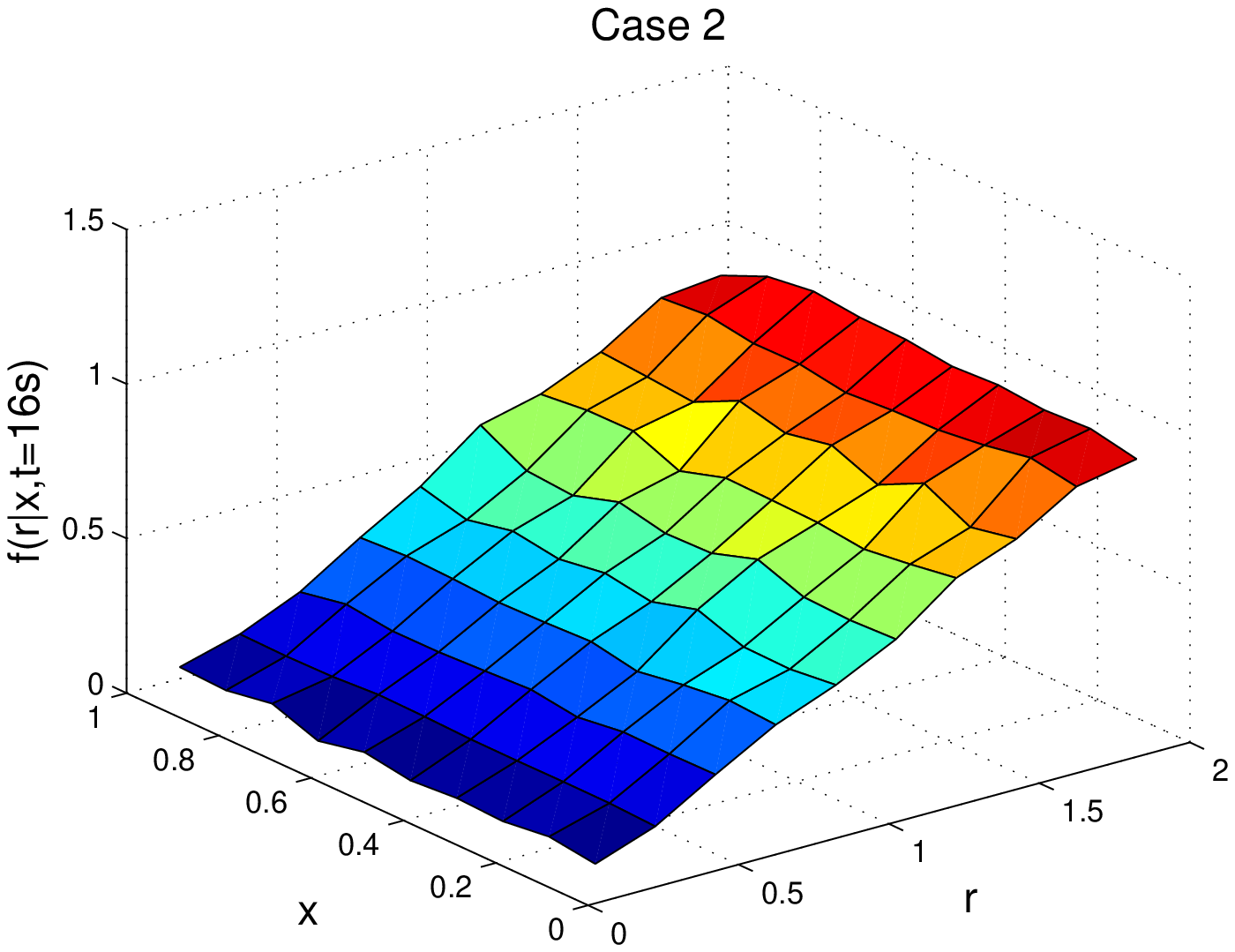,width=\textwidth}
\end{minipage}
\hfill
\begin{minipage}{4.5cm}
\epsfig{figure=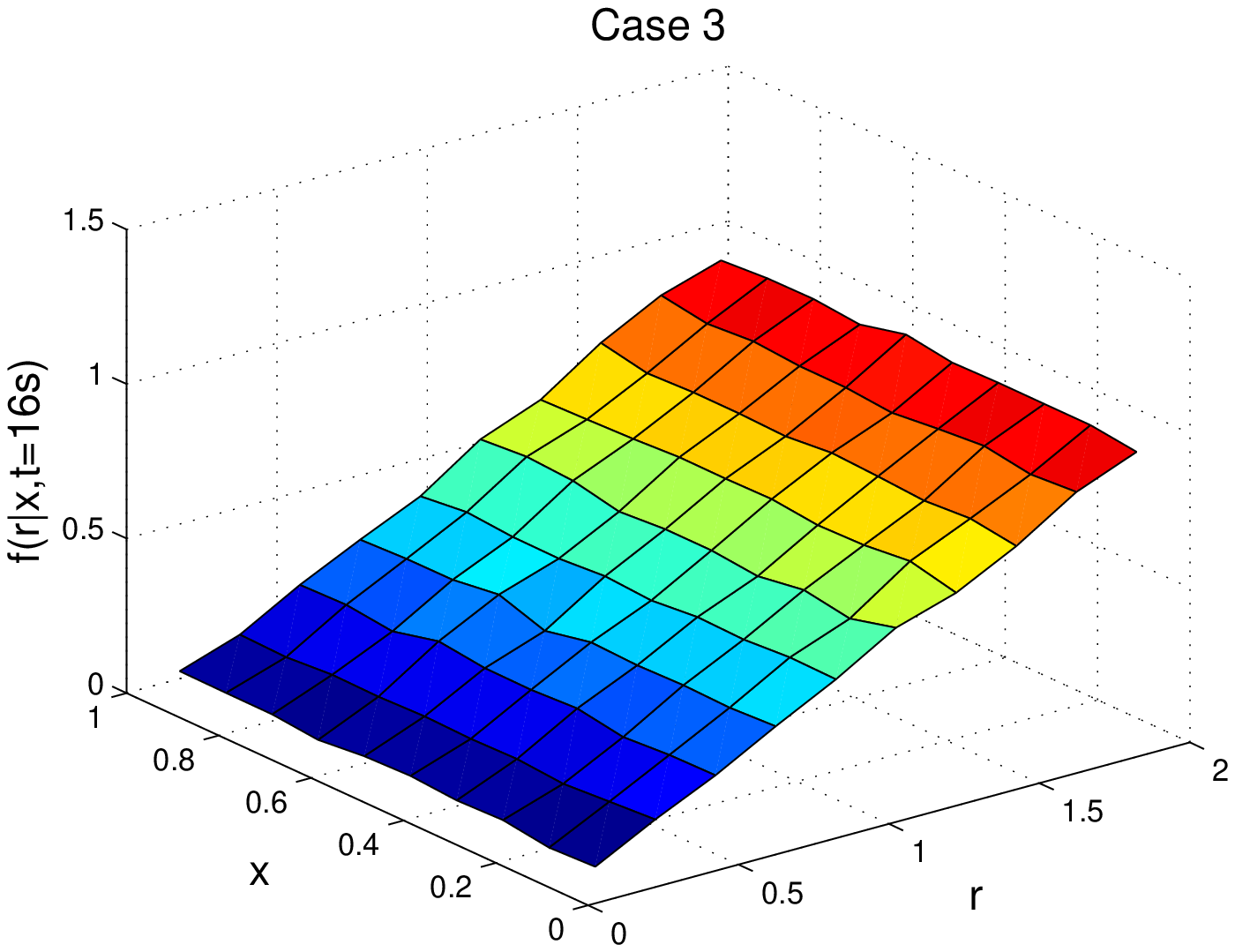,width=\textwidth}
\end{minipage}
\hfill
\begin{minipage}{4.5cm}
\epsfig{figure=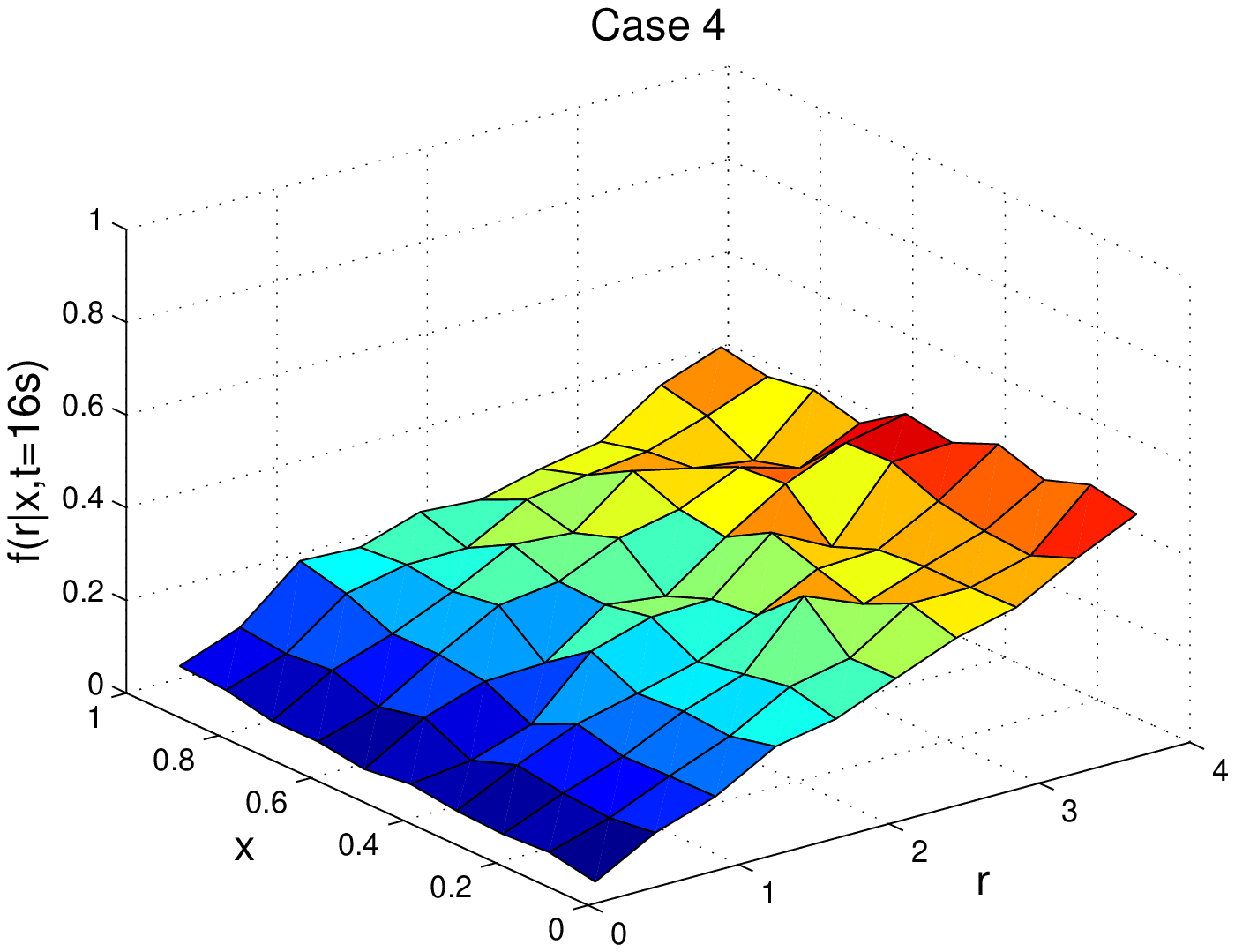,width=\textwidth}
\end{minipage}
\hfill
\begin{minipage}{4.5cm}
\epsfig{figure=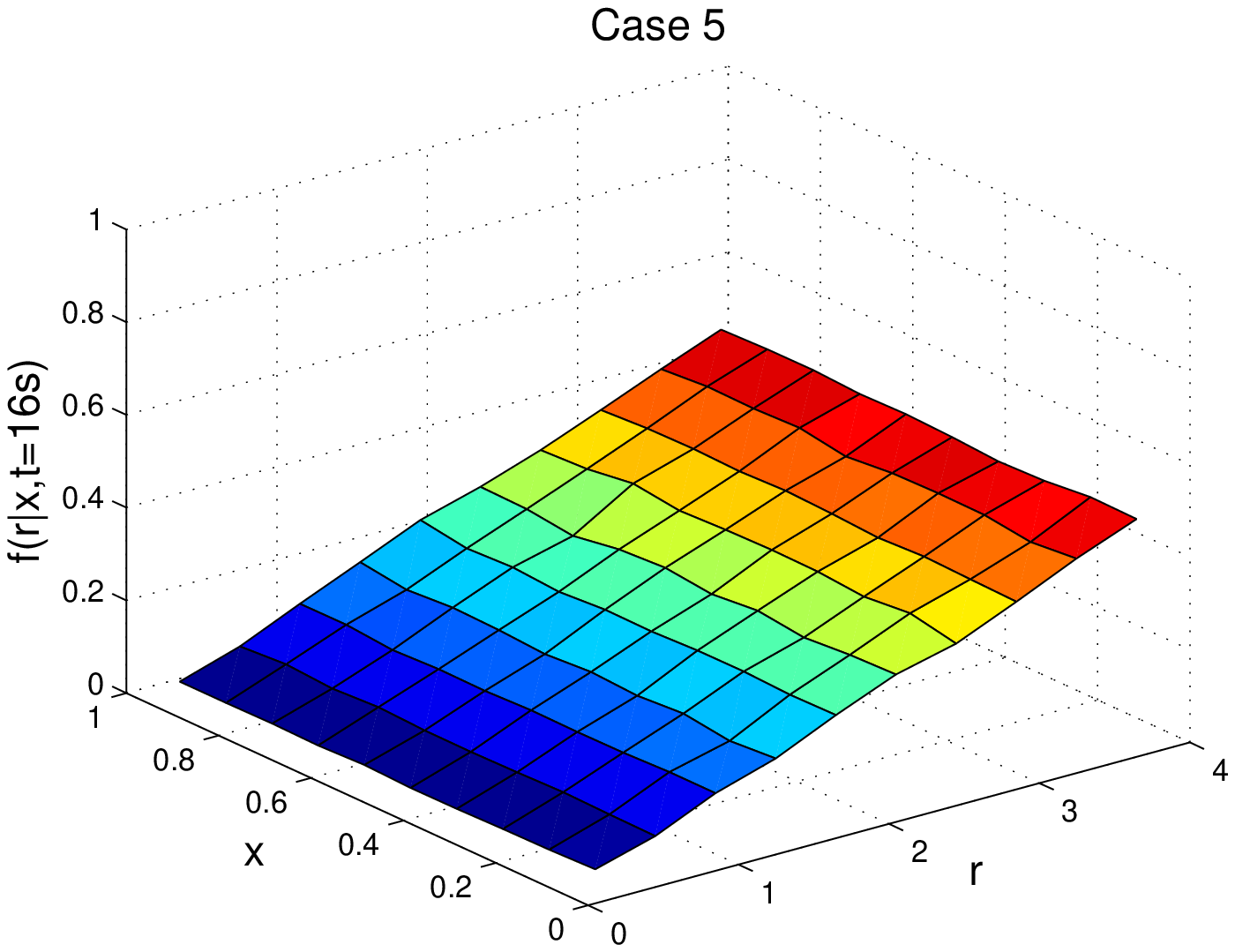,width=\textwidth}
\end{minipage}
\hfill
\begin{minipage}{4.5cm}
\epsfig{figure=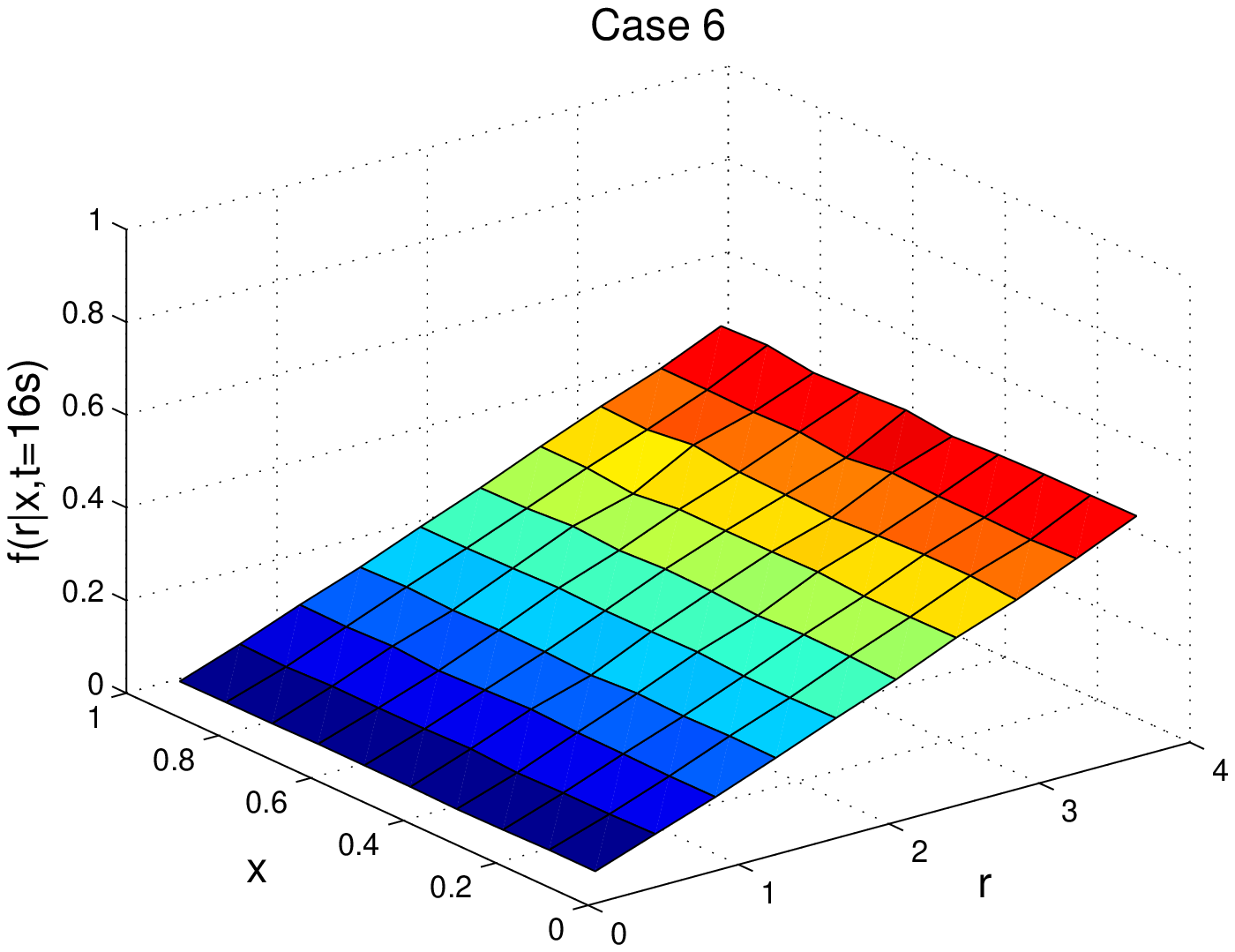,width=\textwidth}
\end{minipage}
\hfill
\begin{minipage}{4.5cm}
\epsfig{figure=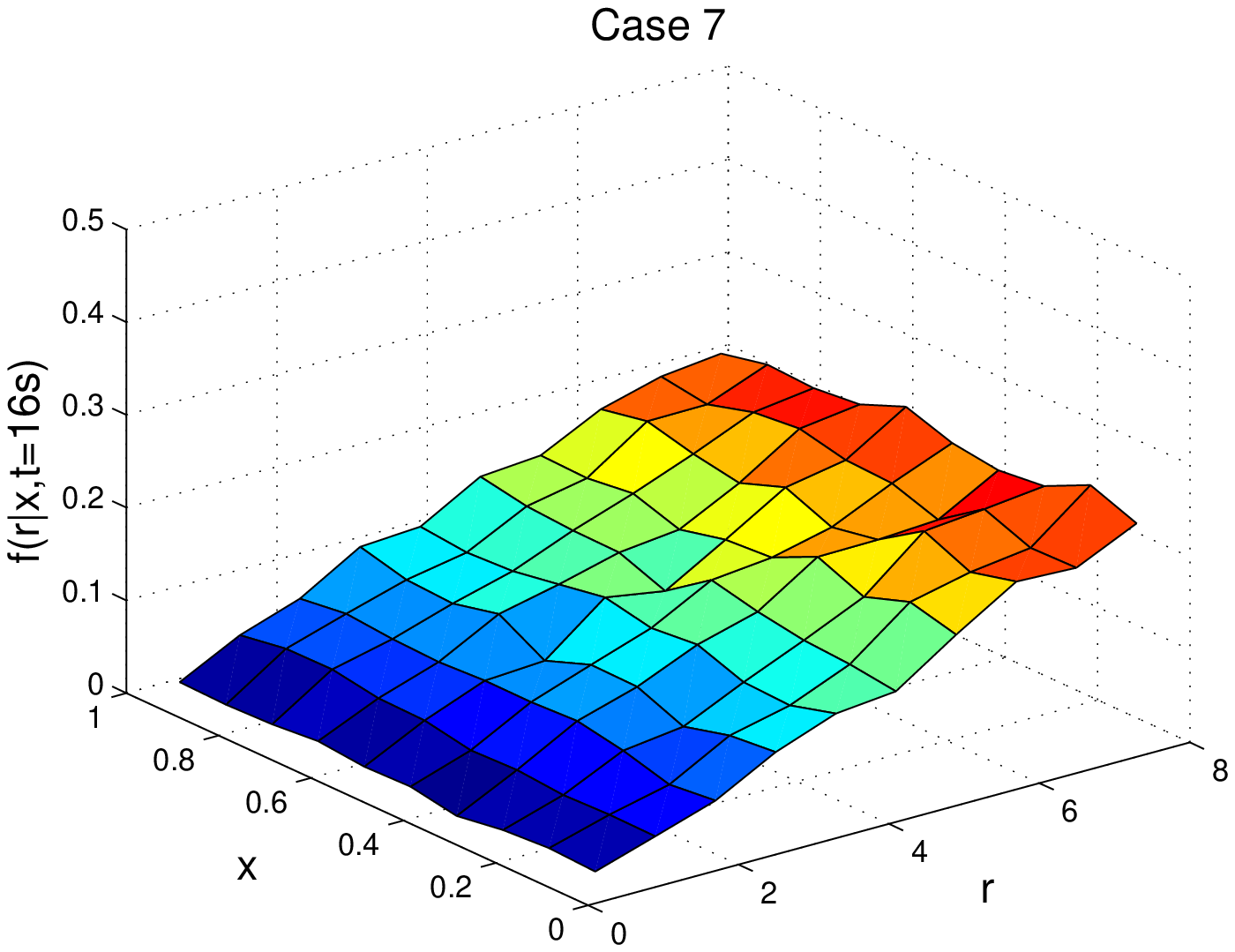,width=\textwidth}
\end{minipage}
\hfill
\begin{minipage}{4.5cm}
\epsfig{figure=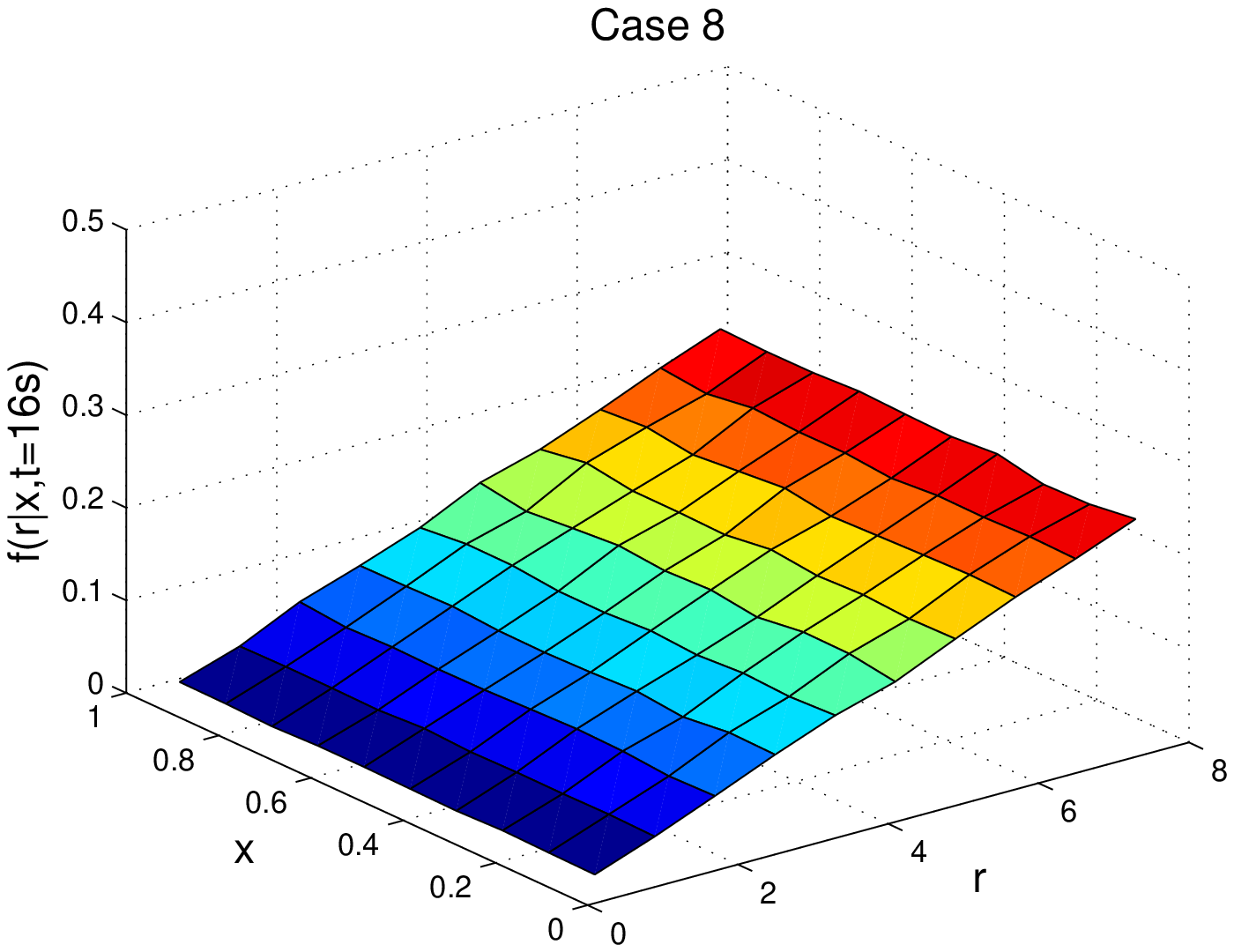,width=\textwidth}
\end{minipage}
\hfill
\begin{minipage}{4.5cm}
\epsfig{figure=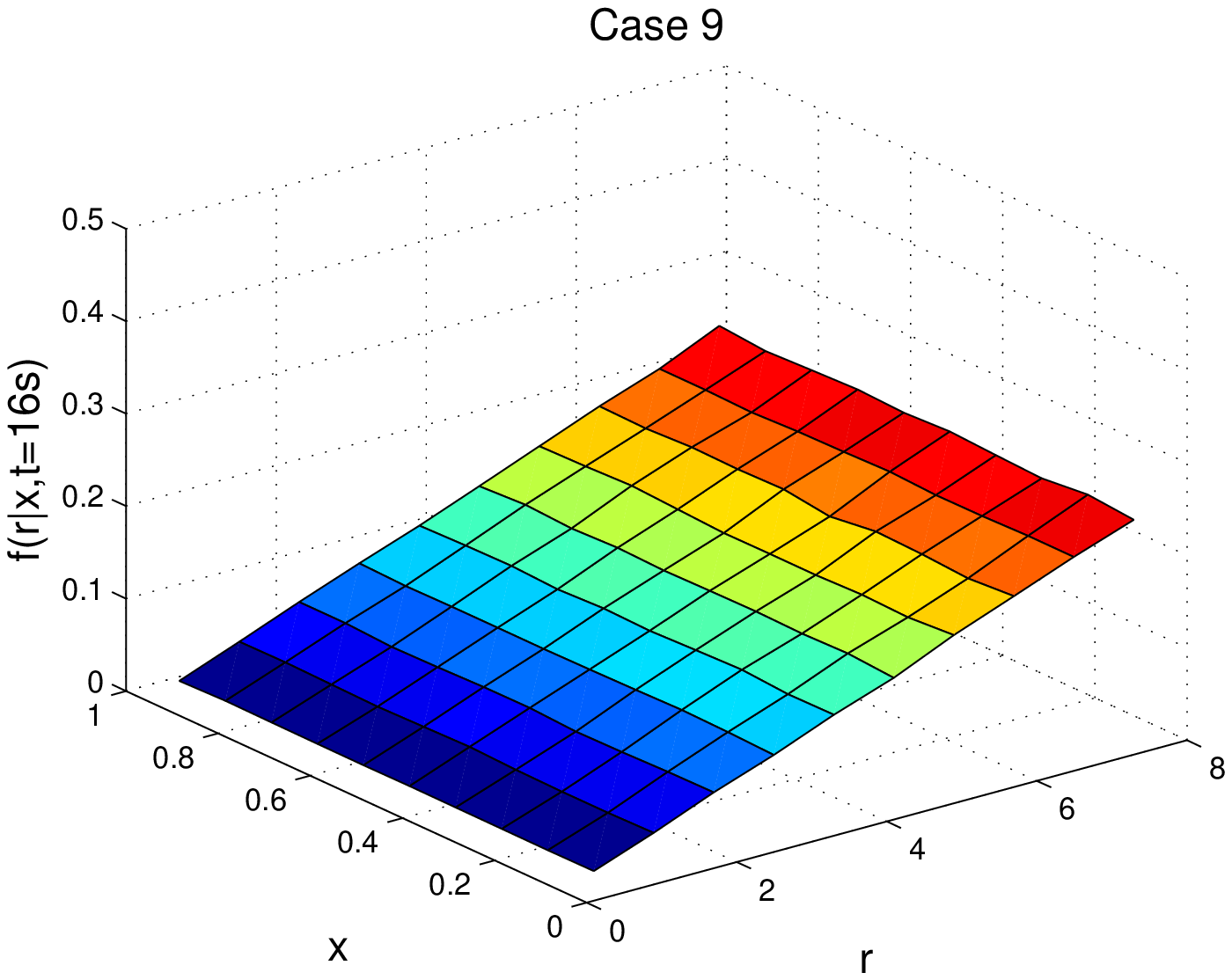,width=\textwidth}
\end{minipage}
\hfill
\caption{Conditional distributions of TPT, $f(r|x,t)$}
\label{condden:fig}
\end{figure}

\begin{figure}
\begin{minipage}{8cm}
\center \epsfig{figure=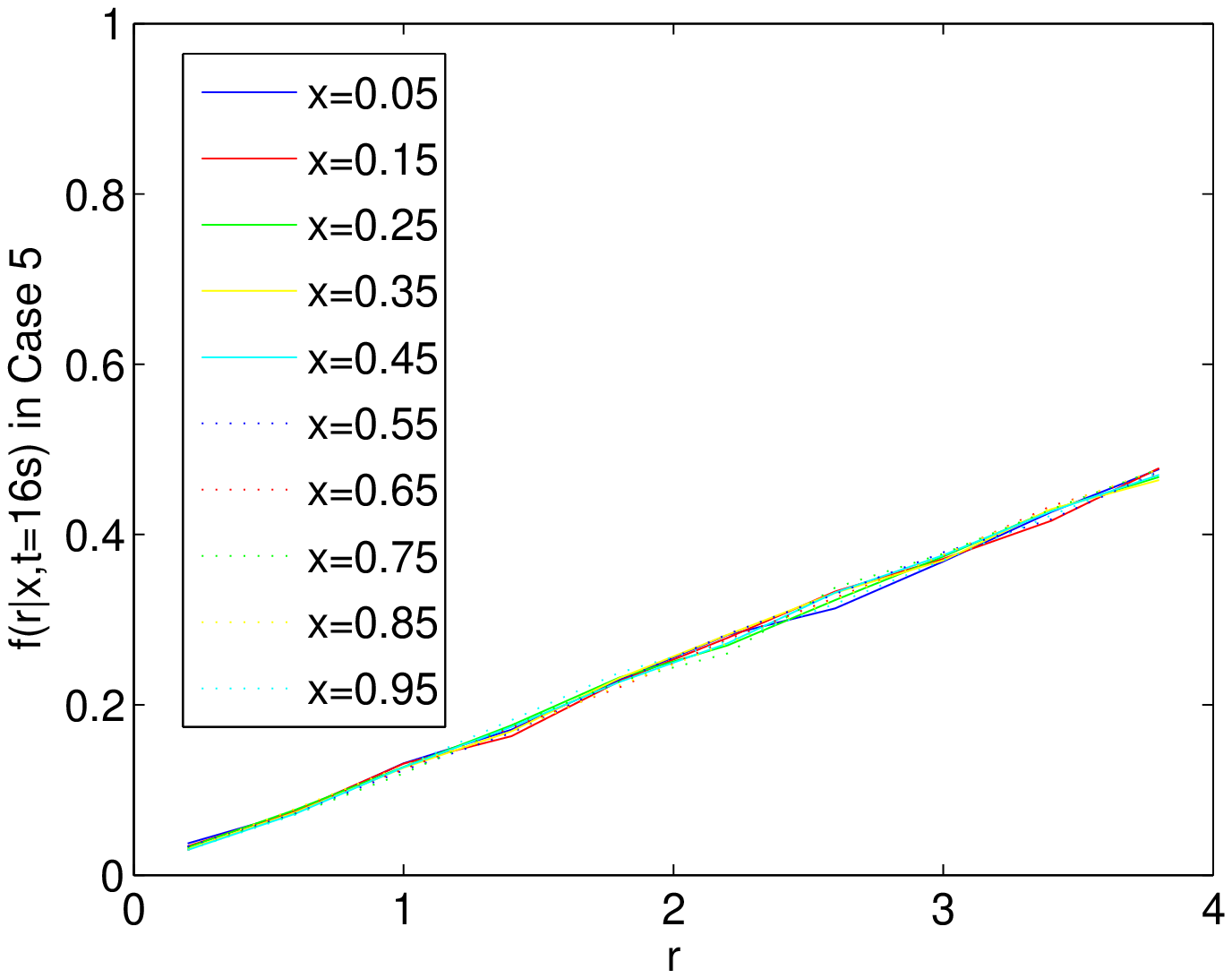,width=\textwidth}
\end{minipage}
\hfill
\begin{minipage}{8cm}
\center \epsfig{figure=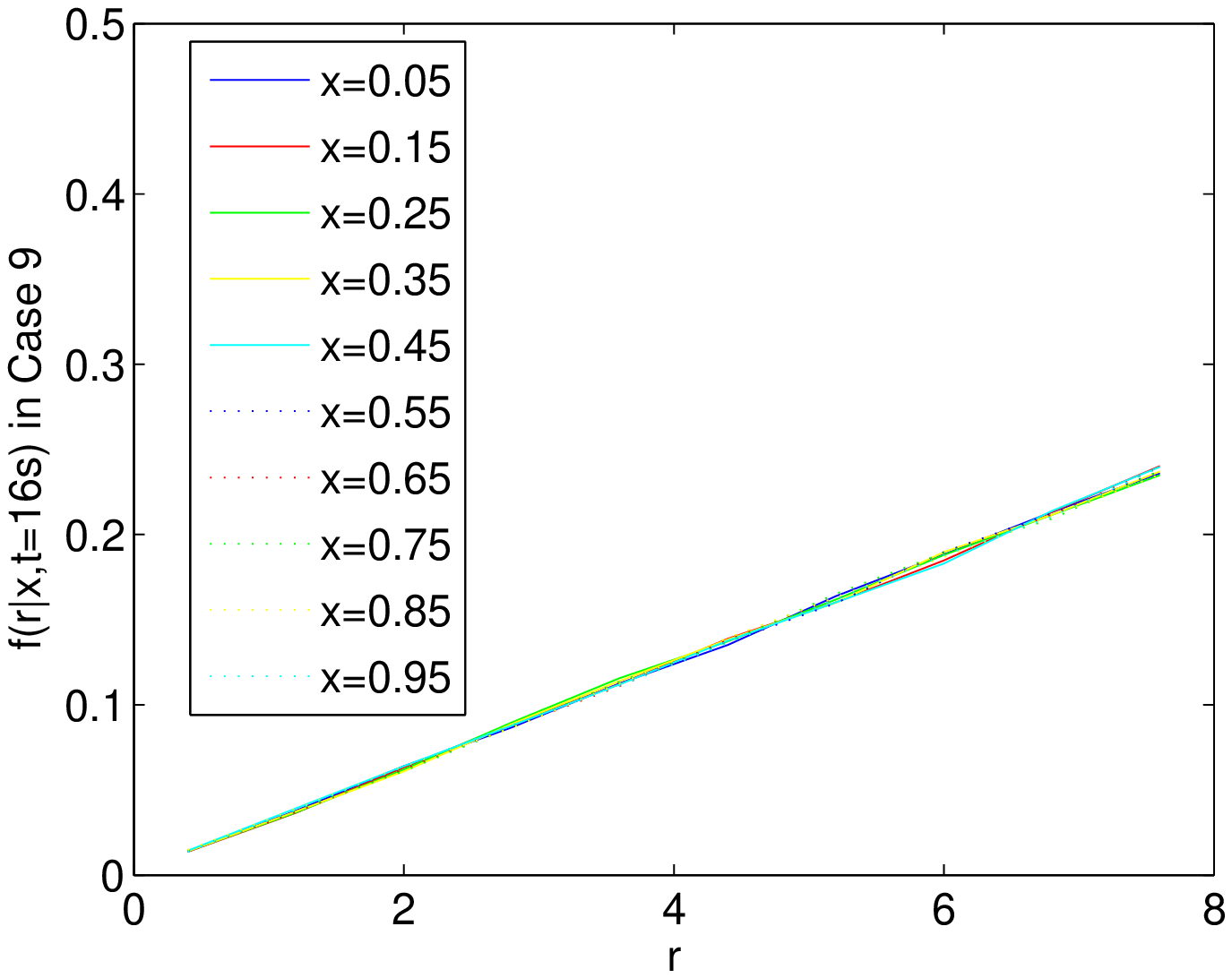,width=\textwidth}
\end{minipage}
\hfill
\caption{Conditional TPT distribution, $f(r|x,t=16s)$. Top: Case 5; Bottom: Case 9}
\label{conddencom:fig}
\end{figure}

In what follows, the discrete model (\ref{discrete2:eqn}) is
executed with constant influxes and uniform TPT distributions given
in Table 1.
For all cases, an ensemble of $7,000$ realizations is used in order
to obtain smooth evolution of the densities $f(x,r,t)$ and
$\rho(x,t)$.
At time t=16{\it sec}, the phase $\phi$ and TPT $\tau$ of items
existing in the factory (i.e., $0 \le \phi <1$) are recorded and
used to plot a surface for the conditional density $f(r|x,t)$ (Fig.
\ref{condden:fig}).
It is found that when $\lambda_0 \sigma_{\mathcal T}^0$ is sufficiently
large (as in Case 3,5,6,8,9), $f(r|x,t)$ has become independent of
the phase coordinate $x$.
Actually, for our particular choice of ${\mathcal T}(r,t)$ shown in cases 3,5,6,8,9 of 
Table 1, we find that the conditional density $f(r|x,t)$ follows a linear
distribution in $r$ passing through the origin (Fig.
\ref{conddencom:fig}).
Clearly, in order to lift effectively, one needs to know the {\it
support} of the conditional distribution; for our choice of
${\mathcal T}(r,t)$ the lower limit of its support in $r$ does not depend
on WIP or time, and it makes sense to take this also as the lower
limit of support of $f(r|x,t)$; we have observed that - to within
acceptable error - the upper limit of support of $f(r|x,t)$ coincided
in our simulations with the upper limit of
support of ${\mathcal T}(r,t)$.
In this case, $f(r|x,t)$ can be approximated in the form of 
$r {\mathcal T} (r,t)/C$ where $C$ is a constant independent of $r$.
Obviously, $C$ equals $T_1(t)$ since $\int f(r|x,t) = 1$. Other choices of $\lambda(t)$ and ${\mathcal T}(r,t)$ (e.g., a linearly
increasing influx and a linear TPT distribution) give rise to more
or less the same relationship between $\rho(x,t)$ and $f(x,r,t)$.
The ``constitutive equation'' (\ref{reducedrelation:eqn}) and the
parameter regime under which it is justified were thus found using
only short simulations with the discrete fine-level model.
\begin{table}
\begin{center}
\begin{tabular*}{0.75\textwidth}%
   {@{\extracolsep{\fill}}cccr}
   Case Number  & influx, $\lambda(t)$  & PDF of TPT, ${\mathcal T}(r,t)$ \cr
   \hline
   1                  &  0.5      & uniform in [0.1,2]  \cr
   2                  &  10       & uniform in [0.1,2]  \cr
   3                  &  20       & uniform in [0.1,2]  \cr
   4                  &  0.5      & uniform in [0.1,4]  \cr
   5                  &  10       & uniform in [0.1,4]  \cr
   6                  &  20       & uniform in [0.1,4]  \cr
   7                  &  0.5      & uniform in [0.1,8]  \cr
   8                  &  10       & uniform in [0.1,8]  \cr
   9                  &  20       & uniform in [0.1,8]  \cr
   \hline
\end{tabular*}
\end{center}
\hfill
\hfill
{Table 1 Cases of influxes and uniform distributions of TPT both constant over the time domain}
\label{influxTPT:tab}
\end{table}

For a prescribed ${\mathcal T}(r,t)$, and starting at a given
$\rho(x,t)$, the lifting algorithm can now be formulated when
$\lambda_0 \sigma_{\mathcal T}^0 \gg 1$ as follows:
\begin{enumerate}
\item Calculate the WIP(t): $WIP = \int_0^1 \rho(x,t) dx$.
\item Calculate the probability density function (or normalized number density)
$f_{\rho}(x,t)$
of the phase coordinate: $f_{\rho}(x,t) = {{\rho(x,t)} \over {WIP}}$.
\item Since the number of items generated can only be an integer,
we have to systematically generate an ensemble of integers whose
mean value equals WIP.
Let $a=int(WIP)+1-WIP$, where $int(WIP)$ is
the maximum integer not greater than WIP.
Select a random variable
$p$ which is uniformly distributed in $[0,1]$.
If $p<a$, then the
number of items is $int(WIP)$.
Otherwise it is $int(WIP)+1$.
\item Compute the cumulative distribution function (CDF) of the
phase coordinate $F_X(x,t)=\int_0^x f_{\rho}(x,t)$ and its inverse $IF_X(F_X,t)$.
\item For each realization in the integer ensemble,
generate the phase coordinates of items:
$\phi_i = IF_X(\zeta_i,t), i=1,2,\cdots,n$,
where $n$ is the integer number of items and $\zeta_i,i=1,2,\cdots,n$
are random numbers uniformly distributed in $[0,1]$.
\item Compute the CDF of the TPT $F_R(r,t)=\int_{-\infty}^r f(r|x,t)$
and its inverse $IF_R(F_R,t)$.
The phase coordinate $x$ is suppressed in $F_R(r,t)$ since it is
independent of $x$.
\item For items in each realization, generate their
respective TPT's: $\tau_i$ $=$ $IF_R(\psi_i,t)$, $i=1,2,\cdots,n$,
where $\psi_i$ are also random numbers uniformly distributed in $[0,1]$.
\end{enumerate}

The following test was performed to validate the above procedure.
A
true trajectory of $\rho(x,t)$ for Case 9 in Table
1 is computed using $5,000$ ensemble realizations
in the time domain directly by the discrete phase model
(\ref{discrete2:eqn}).
Then the trajectory is interrupted at $t=10$;  the number density
$\rho(x,t)$ at the interruption is retained, and then lifted as
described above to start a new discrete evolution.
We observed that the relative difference between the 
two trajectories (the normally
continued one and the one starting from the lifting after the interruption) is quite small (within
$2\%$) even just after the lifting.
Another realization of the lifting may be taken which uses a
different random number seed.
The result shows no significant difference than the previous situation.
It can therefore be concluded that the lifting algorithm is
effective and particular liftings do not significantly affect
subsequent restricted number densities.

{\it 3.2$\quad$Coarse Projective Integration of the Density-level System} \label{EFevolution:sec}

Coarse Projective Integration (CPI) is a numerical technique
developed in the EF framework to evolve in time the coarse-grained
observables of a multiscale system, estimating their temporal
derivatives using the underlying fine-scale simulator
\cite{Gear:01,GearA:01,Rico-Martinez:04}.
This technique is suitable for systems whose coarse- and fine-level
temporal scales are well separated, i.e., the coarse-level
observables are smooth over a temporal scale that is significantly
larger than the fine-level evolution scales (the scales that it
takes for higher-order system observables to become slaved to the
slow, ``master" ones).
In the same spirit with adaptive time-step selection based on error
control for standard deterministic integrators \cite{Press:92},
real-time computational tests for projective step selection have to
be taken to ensure error control in the CPI evolution of macroscopic
observables.

\begin{figure}
\begin{minipage}{8cm}
\epsfig{figure=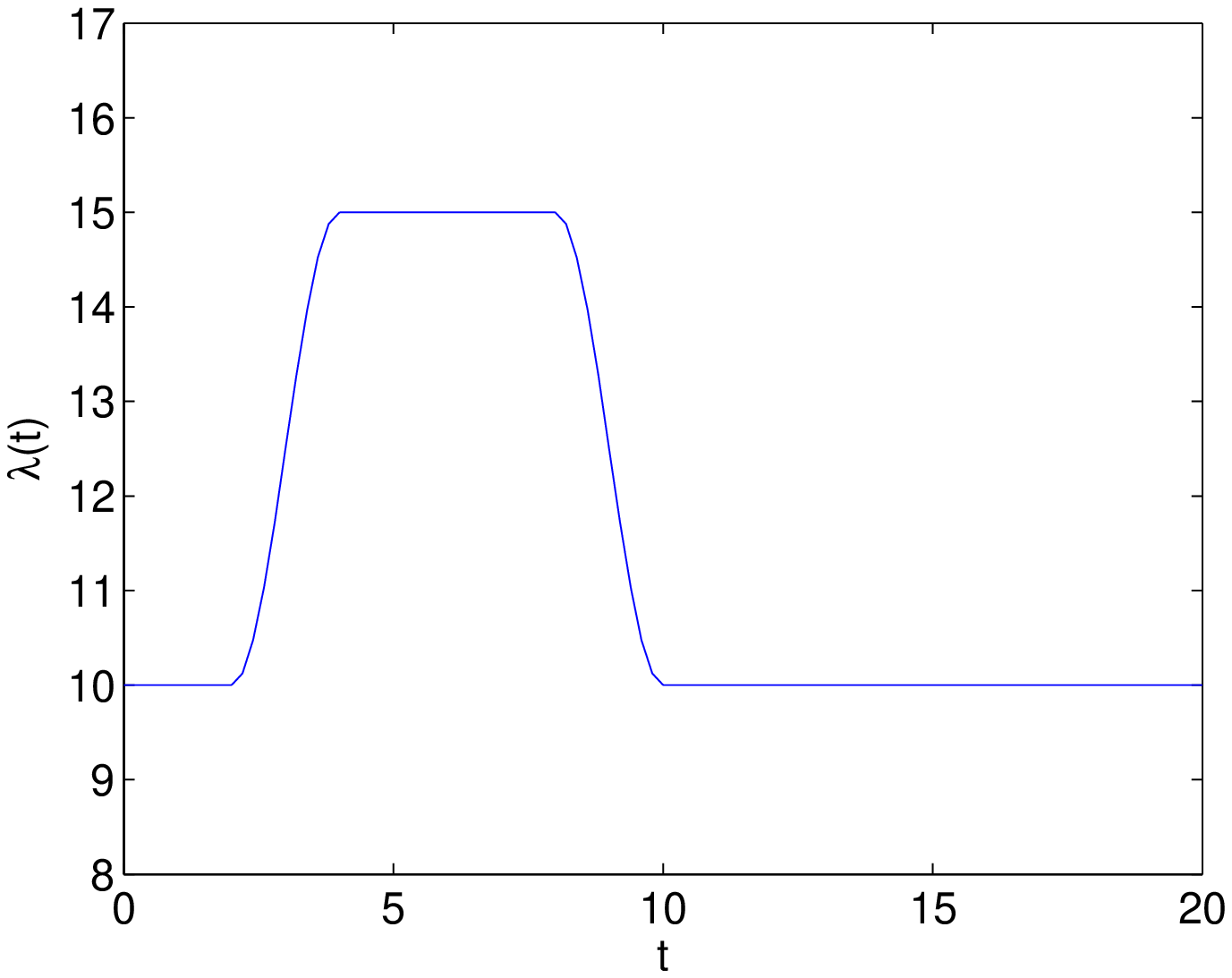,width=\textwidth}
\end{minipage}
\hfill
\begin{minipage}{8cm}
\epsfig{figure=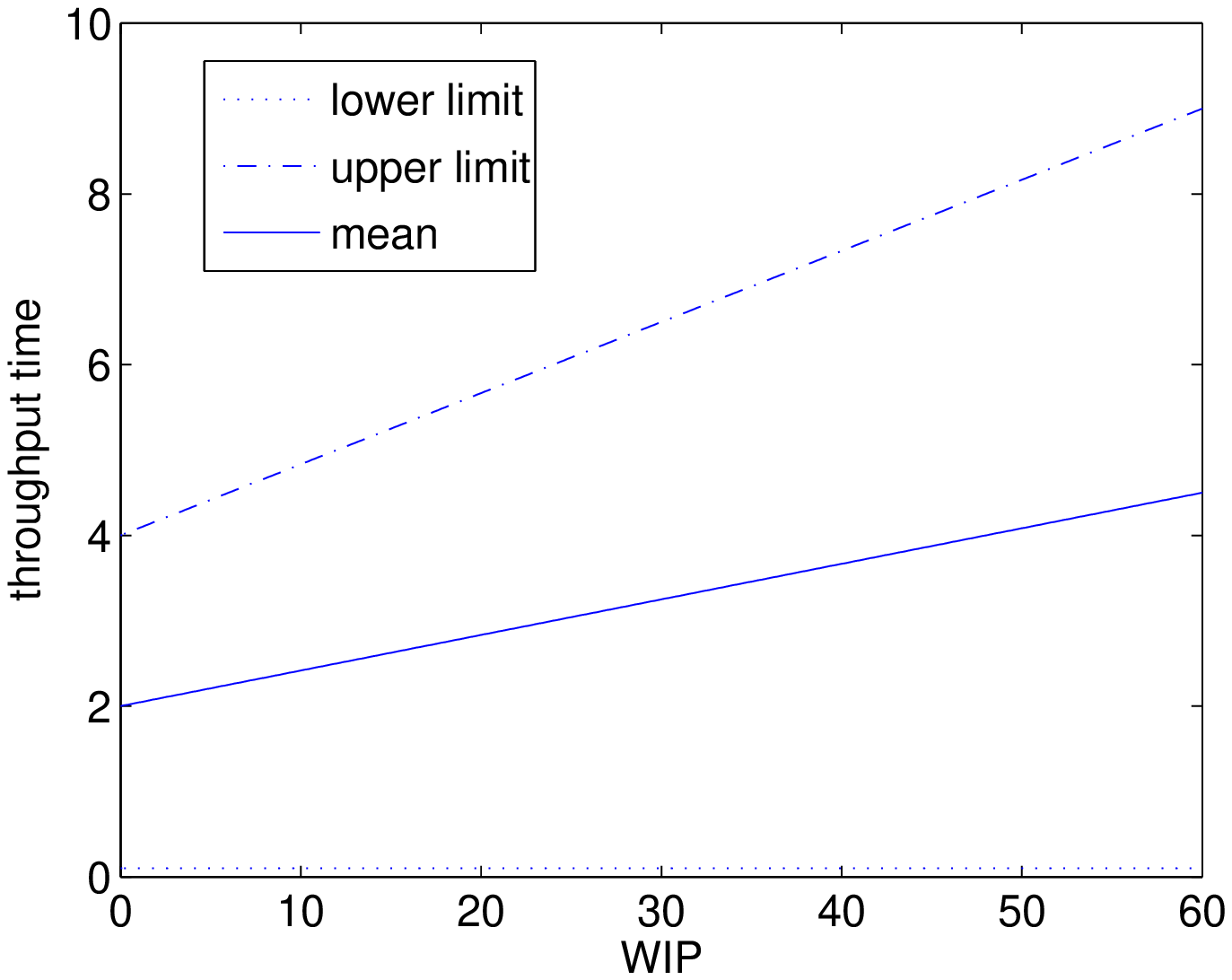,width=\textwidth}
\end{minipage}
\hfill
\caption{Influx and probability distribution of the throughtput time}
\label{CPIinfluxTPT:fig}
\end{figure}

We now implement CPI for our two-scale supply chain model.
The prescribed influx and TPT distribution are shown in Fig.
\ref{CPIinfluxTPT:fig}, respectively.
The time
step $\Delta t_d$ in the discrete model is chosen as $\delta t_d =
10^{-3} sec$.
In the following, a finite-
difference motivated representation is used to evolve coarse-grained
observables, $\rho(x,t)$, through CPI. 
If an explicit equation for $\rho(x,t)$ was discretized in space
with finite differences, we would evolve, in time, the values of the
field at a number of points -- and we would use differences
between these values at each moment in time to approximate the 
{\it right-hand side} of the explicit equation.
Here, we do not use the finite point representation to approximately
evaluate the right-hand side of the equation -- we use the
representation instead to {\it lift} to an item distribution in the
entire domain, and run the fine-scale simulator to {\it estimate}
the {\it left-hand side} of the evolution of the representation.
It is these values that integration codes need to solve the initial
value problem, and the way they process the numbers does not depend
on whether they came from approximating the right-hand side of the
equation, or from estimating its left-hand side.
We projectively integrate the number density $\rho(x,t)$ at equally
spaced phase points $x_j=j/M,j=0,1,\cdots,M$; here we choose $M =
8$. (We have also successfully performed coarse projective integration using 
as variables
the projections of the solution on the first few of its empirical orthogonal global 
basis functions (POD
modes) in a POD-assisted approach; these results are not shown because 
of space
limitations). At the starting time $t_0$, the number density
$\rho(x_j,t_0),j=0,1,\cdots,8$ is 
lifted (according to the lifting algorithm in Section 3.1) 
to initiate the discrete phase
model (\ref{discrete2:eqn}).
This model is subsequently evolved for 20 discrete time steps.
At steps $i, i=12,14,\cdots,20$, phases of items in progress are
restricted and binned to obtain number density histories
$\rho(x,t_0+i\delta t_d)$.
These number density histories are used to approximate the temporal
derivatives $\partial \rho(x,t)/ {\partial t} $ at time
$t_0+20\delta t_d$ via least-squares fitting.
The number density after a coarse-level time interval $\Delta t_c$
can then be obtained based on a simple forward Euler explicit
algorithm by {\it projection}:
$$
   \rho(x,t_0+\Delta t_c) = \rho(x,t_0+20\delta t_d) + (\Delta t_c-20\delta_d) { {\partial \rho(x,t_0+20\delta t_d)} \over {\partial t}}.
$$
This explicit Euler scheme illustrates the simplest way to project
the number density over a coarse time interval, here with
computational savings of roughly $(1-20\delta_d/\Delta_c)$.
More sophisticated projective algorithms, including projection
templated on higher-order continuum integration schemes
\cite{Rico-Martinez:04}, taking advantage of multiple
spectral gaps \cite{Gear:03} and even coarse implicit schemes
\cite{Gear:02}, can be utilized as well.

\begin{figure}
\begin{minipage}{8cm}
\epsfig{figure=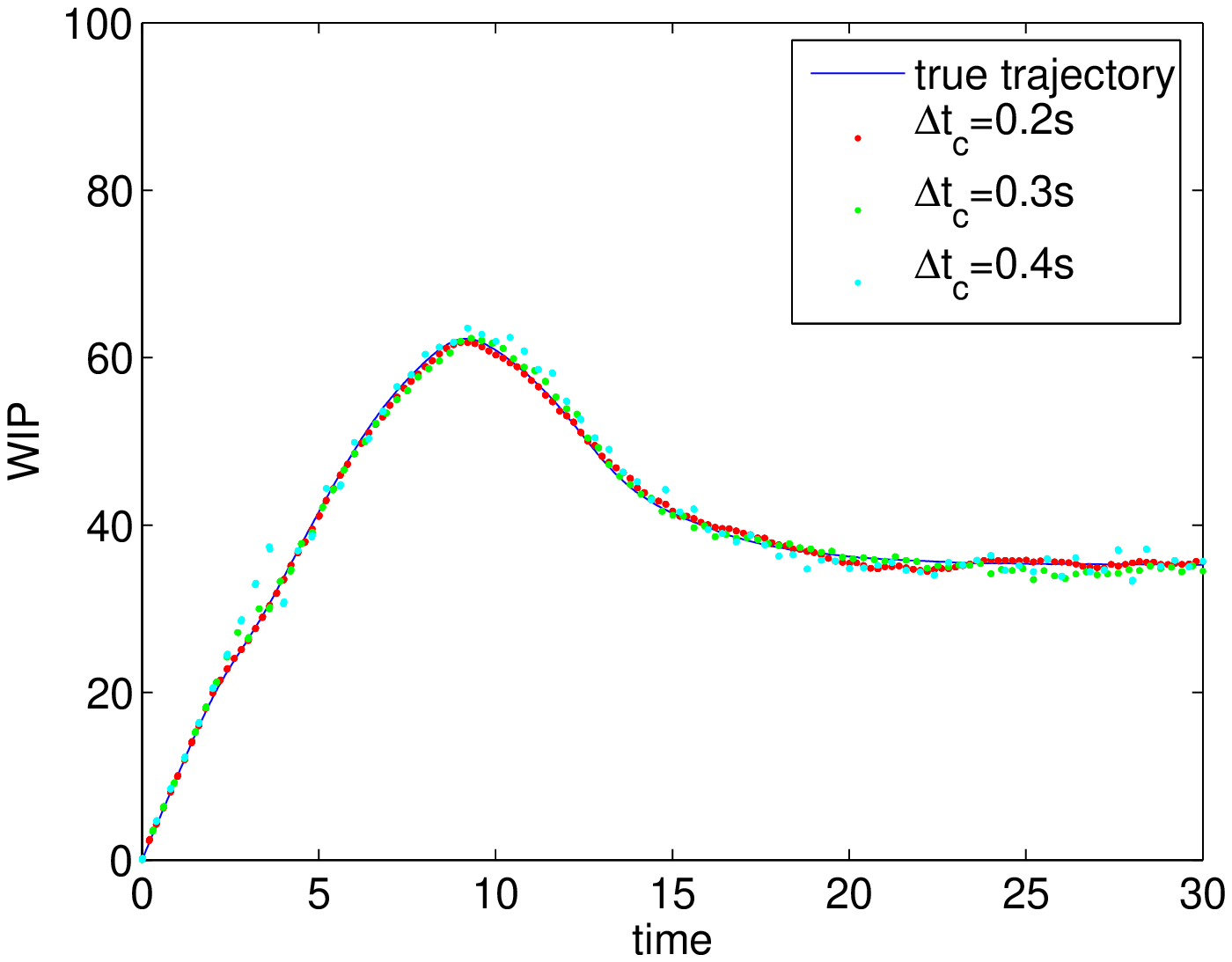,width=\textwidth}
\end{minipage}
\hfill
\begin{minipage}{8cm}
\centerline{\epsfig{figure=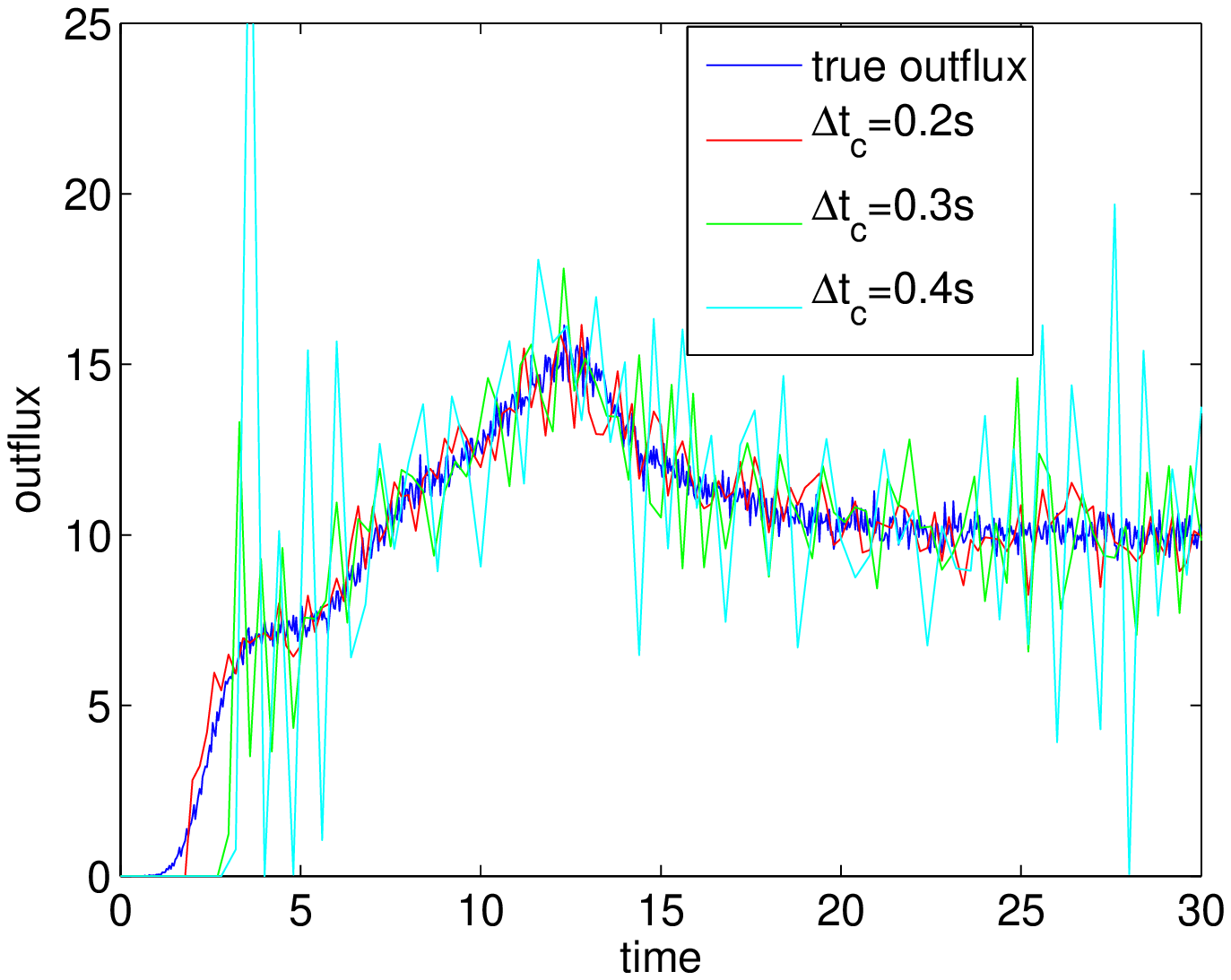,width=\textwidth}}
\end{minipage}
\caption{True trajectories and some CPI evolutions of WIP and ouflux. Left: WIP, clusters of points represent WIP's at time steps $ik\Delta t_d, i=0,1,\cdots,10$ immediately following the lifting; Right: outflux.}
\label{CPIWIPoutflux:fig}
\end{figure}

\begin{figure}
\begin{minipage}{8cm}
\epsfig{figure=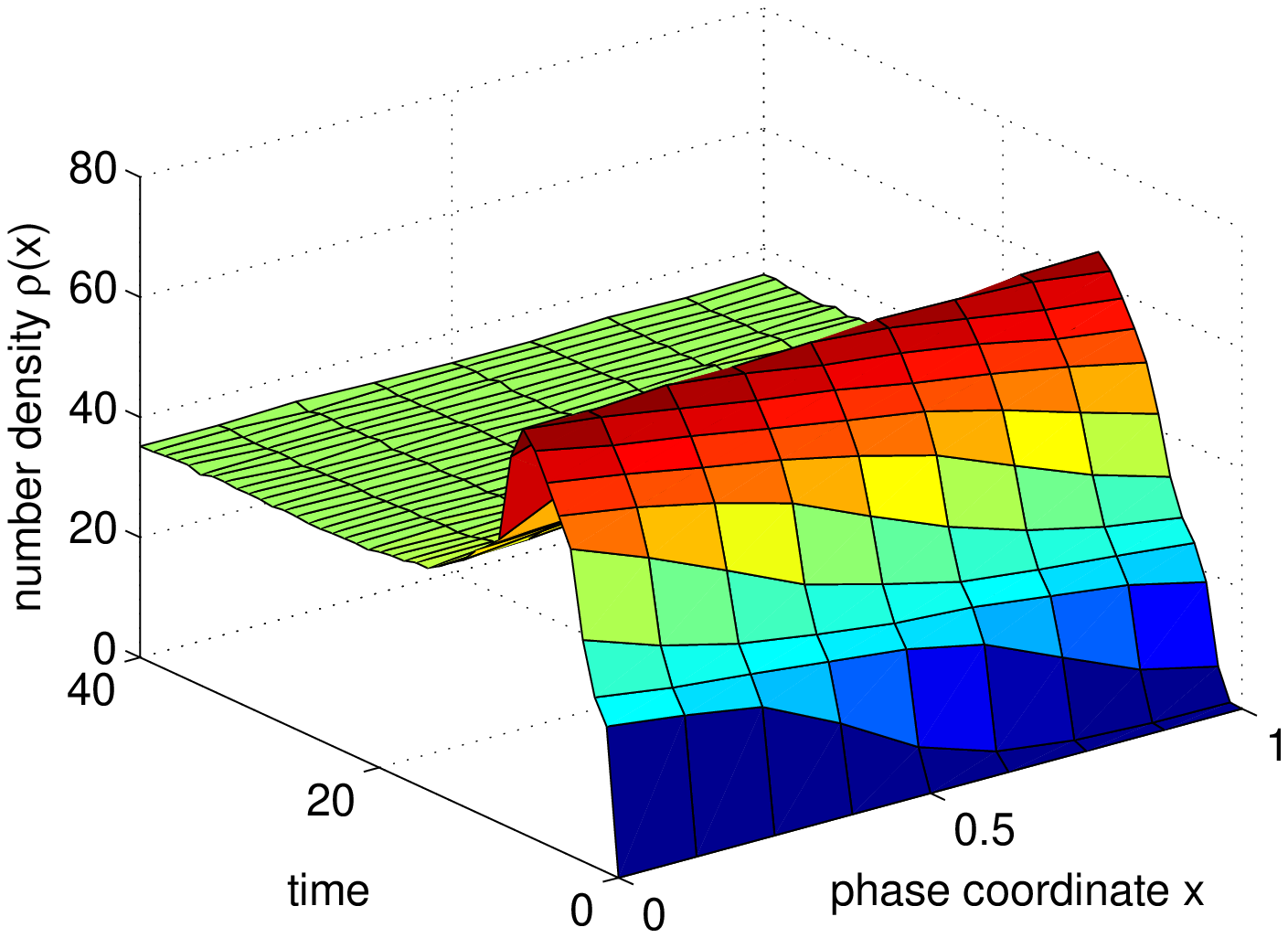,width=\textwidth}
\end{minipage}
\hfill
\begin{minipage}{8cm}
\epsfig{figure=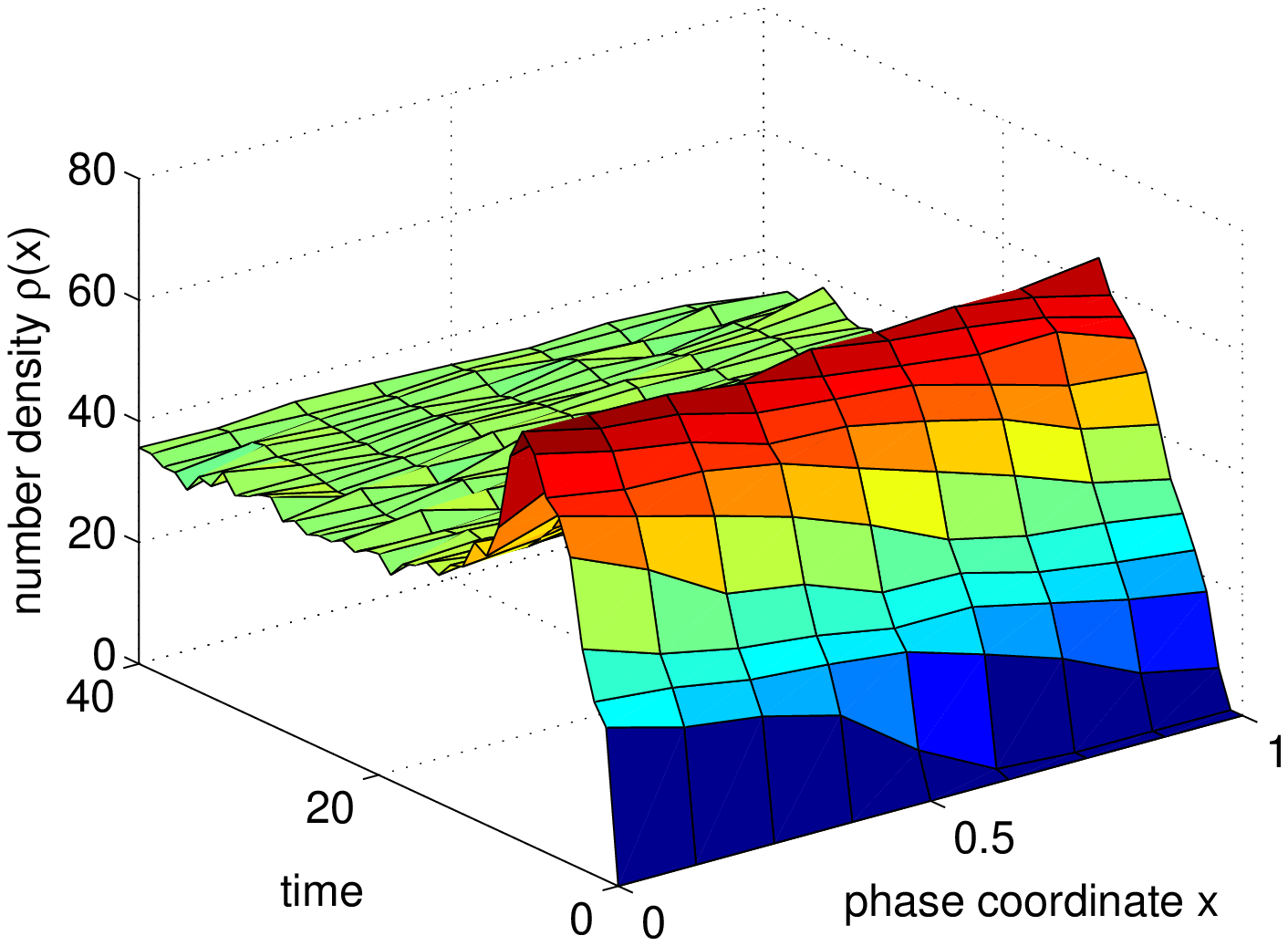,width=\textwidth}
\end{minipage}
\caption{Number densities of the true evolution and CPI evolution with $\Delta t_c=0.2s$. Left: true evolution; Right: CPI result}
\label{CPIrho:fig}
\end{figure}

Figure \ref{CPIWIPoutflux:fig} compares trajectories of WIP and outflux in true and CPI evolutions 
for some selections of the coarse time step $\Delta t_c$. All computations involve an ensemble of $5,000$ realizations. 
It is found that the CPI trajectories of WIP almost coincide as $\Delta t_c = 0.2s$ and $\Delta t_c = 0.3s$, which suggests 
that the result at $\Delta t_c = 0.2s$ is correct. Comparison also displays that 
WIP and outflux at $\Delta t_c = 0.2s $ indeed match their true trajectories. Number densities of the true evolution and 
the CPI result with $\Delta t_c=0.2s$ are also shown in Fig. \ref{CPIrho:fig}. By using coarse projective integration, 
the computational load is substantially reduced by $78.3\%$, taking into account the trial efforts in selecting an appropriate coarse 
timestep size for projection.

{\bf 4$\quad$Conclusions} \label{conclusion:sec}

We demonstrated certain features of Equation-Free coarse-grained
computation for a re-entrant supply-chain model.
In cases where an explicit equation for the processed item number
density cannot be easily derived, simulations with the fine-scale
model can be used to determine the conditions under which such an
equation can, in principle, exist. 
Once these conditions are established, equation-free methods such as
coarse projective integration can be used to evolve the coarse-level
density through short bursts of appropriately initialized runs with
the fine-scale simulator.
The lifting process essential to this EF approach was investigated
and a particular governing form establishing a relation between the
joint number density $f(x,r,t)$ and the number density $\rho(x,t)$
was found; the form is of course only valid for this model.

EF algorithms based on matrix-free iterative linear algebra (like
Newton-Krylov GMRES, \cite{Kelley:95}) can also be used to effectively
implement contraction mappings to find stationary item number
densities, if the influx to a factory reaches a stationary state at long times. Additional
tasks like continuation, stability and parametric sensitivity
analysis, and even control and optimization computations can in
principle be implemented in an equation-free framework.

Our illustrative example involved a Monte-Carlo type solution of a 
Boltzmann equation, which, at a certain limit, approaches a discrete
event simulator.
In this case, both the Boltzmann equation and the reduced equation
for the item number density at the appropriate limit were
analytically available; this allowed us to validate our
equation-free computations.
The real challenge lies in wrapping equation-free algorithms around
true discrete event simulations for realistic processing factory
configurations.
While the derivation of continuum-level equations may be extremely
difficult or practically infeasible in such cases, our
computer-assisted approach remains essentially the same, independent
of the details in the underlying fine-scale simulator.
We believe that this approach holds promise for facilitating the
extraction of system-level information from complex discrete event
simulators.

{\bf 5$\quad$Acknowledgements} 

The research of DA was supported by NSF grant DMS-0204543. 
IGK and YZ gratefully acknowledge support by DOE and by an NSF/ITR grant.







\bibliographystyle{plain}
\bibliography{yannis_paper_6}

\end{document}